\documentclass[a4paper,12pt,reqno]{amsart}

\usepackage[text={130mm,200mm}]{geometry}

\theoremstyle{plain}
\newtheorem{theorem}{Theorem}
\newtheorem*{theorem*}{Theorem}
\newtheorem{corollary}{Corollary}
\newtheorem*{corollary*}{Corollary}
\newtheorem{lemma}{Lemma}
\newtheorem*{lemma*}{Lemma}
\newtheorem{proposition}{Proposition}
\newtheorem*{proposition*}{Proposition}

\newtheorem*{conjecture*}{Conjecture}
\theoremstyle{definition}
\newtheorem{definition}{Definition}
\newtheorem*{definition*}{Definition}
\theoremstyle{remark}
\newtheorem{remark}{Remark}
\newtheorem*{remark*}{Remark}
\newtheorem{example}{Example}

\providecommand{\surname}[1]{#1}
\providecommand{\country}[1]{#1}

\begin{document}

\title[(p-adic)  zeta values]{($p$-adic) $L$-functions and ($p$-adic) (multiple) zeta values}

\author{Nikolaj \surname{Glazunov}}
\address{Department of Electronics\\
  National Aviation University\\
  1~Komarova Pr.\\
  Kiev\\
 03680 \\
  \country{Ukraine}}
\email{glanm@yahoo.com}
\urladdr{https://sites.google.com/site/glazunovnm/}






\begin{abstract}
The article is dedicated to the memory of George Voronoi.
It  is concerned with ($p$-adic) $L$-functions (in partially  ($p$-adic) zeta functions)  and cyclotomic  ($p$-adic) (multiple) zeta values.
The beginning of the article contains a short  summary of the results on the Bernoulli numbers associated with the studies of George Voronoi.
 Results on multiple zeta values have presented by D. Zagier, by P. Deligne and A.Goncharov, by A. Goncharov, by F. Brown, by C. Glanois and others. S. \"Unver have investigated p-adic multiple zeta values in the depth two. Tannakian interpretation of p-adic multiple zeta values is given by H.  Furusho. 
Short history and connections among Galois groups, fundamental groups, motives and arithmetic functions are presented in the talk by Y. Ihara.
Results  on multiple zeta values, Galois groups and geometry of modular varieties has presented by Goncharov.
Interesting unipotent motivic fundamental group is defined and investigated by Deligne and Goncharov.
The framework of ($p$-adic) $L$-functions and ($p$-adic) (multiple) zeta values is based on Kubota-Leopoldt $p$-adic $L$-functions and arithmetic $p$-adic $L$-functions by Iwasawa. Motives and  ($p$-adic) (multiple) zeta values by Glanois and by \"Unver, improper intersections of Kudla-Rapoport divisors and Eisenstein series by Sankaran are reviewed.
More fully the content of the article can be found at the following table of contents: Introduction. 1. Voronoi-type congruences for  Bernoulli numbers. 2. Riemann zeta values. 3. On class groups of rings with divisor theory.  Imaginary quadratic and cyclotomic fields. 4. Eisenstein Series. 5. Class group, class fields and zeta functions. 6. Multiple zeta values. 7. Elements of non-Archimedean local fields and $ p-$adic analysis. 8. $ p-$adic  interpolation of zeta and $L$-functions. 9. Iterated integrals and (multiple) zeta values. 10. Formal groups and $p$-divisible groups. 11. Motives and ($p$-adic) (multiple) zeta values. 12. On the Eisenstein series associated with Shimura varieties.
Sections 1-10 and subsection 12.1 (On some Shimura varieties and Siegel modular forms) can be considered as an elementary introduction to the results of section 11 and subsection 12.2 (On improper intersections of Kudla-Rapoport divisors and Eisenstein series).
Numerical examples are included.
\end{abstract}

\maketitle



\newpage
\section*{Introduction}

The article is dedicated to the memory of George Voronoi.
It  is concerned with ($p$-adic) $L$-functions (in partially  ($p$-adic) zeta functions)  and cyclotomic  ($p$-adic) (multiple) zeta values.
The beginning of the article contains a short  summary of the results on the Bernoulli numbers associated with the studies of George Voronoi.
 Results on multiple zeta values have presented by D. Zagier~\cite{Zagier1992}, by P. Deligne and A.Goncharov~\cite{DeligneGoncharov2005}, by A. Goncharov~\cite{Goncharov2000}, by F. Brown~\cite{Brown2012}, 
 by C. Glanois~\cite{Glanois2016} and others.
Tannakian interpretation of $p$-adic multiple zeta values is given by H. Furusho~\cite{Furusho2007}.
Short history and connections among Galois groups, fundamental groups, motives and arithmetic functions are presented in the talk by Y. Ihara~\cite{Ihara1990}.
Results  on multiple zeta values, Galois groups and geometry of modular varieties has presented by 
Goncharov~\cite{Goncharov2000}.
Interesting unipotent motivic fundamental group is defined and investigated by Deligne and Goncharov~\cite{DeligneGoncharov2005}.
S. \"Unver~\cite{Unver2004,Unver2016}  have investigated  $p$-adic multiple zeta values in the depth two.
The framework of ($p$-adic) $L$-functions and ($p$-adic) (multiple) zeta values is based on Kubota-Leopoldt $p$-adic $L$-functions~\cite{KubotaLeopoldt1964} and arithmetic $p$-adic $L$-functions by Iwasawa~\cite{Iwasawa2001}.
Motives and  ($p$-adic) (multiple) zeta values, improper intersections of Kudla-Rapoport divisors and Eisenstein series by 
Sankaran~\cite{Sankaran2017} are reviewed.
More fully the content of the article can be found at the following table of contents: Introduction. 1. Voronoi-type congruences for  Bernoulli numbers. 2. Riemann zeta values. 3. On class groups of rings with divisor theory.  Imaginary quadratic and cyclotomic fields. 4. Eisenstein Series. 5. Class group, class fields and zeta functions. 6. Multiple zeta values. 7. Elements of non-Archimedean local fields and $ p-$adic analysis. 8. $ p-$adic  interpolation of zeta and $L$-functions. 9. Iterated integrals and (multiple) zeta values. 10. Formal groups and $p$-divisible groups. 11. Motives and ($p$-adic) (multiple) zeta values. 12. On the Eisenstein series associated with Shimura varieties.
Sections 1-10 and subsection 12.1 (On some Shimura varieties and Siegel modular forms) can be considered as an elementary introduction to the results of section 11 and subsection 12.2 (On improper intersections of Kudla-Rapoport divisors and Eisenstein series).
Numerical examples are included.

The subject matter of this review has deep historical roots, with contributions of many mathematiciens. I apologize  for any oversights and any misrepresentations, which are not intentional but rather due to my ignorance. 
\begin{remark*}
Let me now present very briefly the background of my interest on the subject of the values of zeta and $L-$functions.
In 1970-1971 years Yu. Manin gave courses of lectures and seminars on Algebraic Geometry, Diophantine Geometry in MGU and in Steklov mathematical institute. In his lectures and talks Yu. Manin presented and discussed the Birch-Swinnerton-Dyer conjecture concerning $L-$ functions of elliptic curves and abelian varieties. In particular Yu. Manin  have proposed in these talks  modular symbols for computation of values of $L-$functions of elliptic curves at $s = 1$~\cite{Manin1971,Manin1972}. Author of the text attended the lectures and seminars of Yu. Manin. Following of the kind conversation with  Yu. Manin the author has implemented the computer program and has computed Manin's modular symbols~\cite{Glazunov1971} for elliptic curve $E_{\Gamma_0(11)}$ follow to Manin article~\cite{Manin1971}.
\end{remark*}

\section*{Acknowledgements}
I am grateful to J. Steuding for some helpful remarks and comments.

\newpage
\section*{Plan}

{\bf Introduction}
\par\smallskip
{\bf 1. Voronoi-type congruences for  Bernoulli numbers}
\par\smallskip
{\bf 2. Riemann zeta values }
\par\smallskip
{\bf  3. On class groups of rings with divisor theory on imaginary quadratic and cyclotomic fields}
\par\smallskip
{\bf  4. Eisenstein series}
\par\smallskip
{\bf 5.  Class groups, class fields and zeta functions}
\par\smallskip
{\bf 6.  Multiple zeta  values}
\par\smallskip
{\bf 7. Elements of   non-Archimedean local fields and $ p-$adic  analysis}
\par\smallskip
{\bf 8. $ p-$adic  interpolation of zeta and $L$-functions}
\par\smallskip
{\bf 9. Iterated integrals and (multiple) zeta values}
\par\smallskip
{\bf 10. Formal groups and $p$-divisible groups}
\par\smallskip
{\bf  11. Motives and  ($p$-adic) (multiple) zeta values}
\par\smallskip
{\bf 12. On the Eisenstein series associated with Shimura varieties}

\newpage
\section{Voronoi-type congruences for  Bernoulli numbers and generalizations}
We follow to \cite{BorevichShafarevich1986,Voronoi}.
\subsection{Bernoulli numbers} Bernoulli numbers $B_m$ are determined for integers $m \ge 0$ by the expansion
\begin{equation*}
\frac{t}{\exp(t)-1} = 1 + \sum_{m=1}^\infty \frac{B_m}{m!} t^m.
\end{equation*}

\begin{remark}
For $m > 0, B_{2m + 1} =0$.
\end{remark}
So we have $B_0 = 1, B_1 = -\frac{1}{2}, B_2 = \frac{1}{6}, B_4 = -\frac{1}{30}, B_6 = \frac{1}{42}, \ldots$

]smallskip
\subsection{Voronoi's congruences}
Let $N$ be an natural number (the modulus), $a$ coprime with $N$ and let $B_{2m} = \frac{P_{2m}}{Q_{2m}}$ be the Bernoulli number with coprime  $P_{2m}$ and $Q_{2m}$. Then
\begin{equation*}
(a^{2m} - 1)P_{2m} \equiv 2ma^{2m-1} Q_{2m} \sum_{s=1}^{N-1} s^{2m-1} \left[\frac{sa}{N}\right] \bmod N.
\end{equation*}

\subsection{Kummer congruences}
If $p$ is prime and $p-1$ not divide even positive $m$ then the number $\frac{B_m}{m}$ is $p$-integer and there is the congruence
\begin{equation*}
 \frac{B_{m+p-1}}{m+p-1} \equiv \frac{B_{m}}{m} \bmod p.
 \end{equation*}

\subsection{Generalized Bernoulli numbers}
Let $m$ be a natural number.
Let $\chi$ be primitive numeric character (Dirichlet character) modulo $m.$
Generalized Bernoulli numbers $B_{n,\chi}$ are determined for integers $m > 0$ and primitive  Dirichlet character modulo $m$ by the expansion
\begin{equation}
\sum_{a=1}^m \frac{\chi(a)t\exp(at)}{\exp(mt)-1} = 1 + \sum_{n=1}^\infty \frac{B_{n,\chi}}{n!} t^n, \; 0 \le |t| \le \frac{2\pi}{m}.
\end{equation}
\begin{remark}
All the generalized Bernoulli numbers associated with the numerical character $\chi$ are contained in the extension of the field of rational numbers $\mathbb Q$ obtained by joining to $\mathbb Q$ all the values of the character $\chi$.
\end{remark}

\section{Riemann zeta values }

Here we follow to~\cite{Weber1908,Hasse1931,Deuring,BorevichShafarevich1986}.

Let 
$
s = \sigma  + it$ be a complex number and let $\zeta(s) $ be the Riemann zeta function which is presented for $ \sigma > 1$ by the series
\begin{equation*}
\zeta(s) = \sum_{n=1}^\infty \frac{1}{n^s}  
\end{equation*}
;

 By Euler for 
$m \ge 1$   
\begin{equation*}
\zeta(2m) = (-1)^{m-1} \frac{(2\pi )^{2m}}{2(2m)!}B_{2m}
\end{equation*}
 where $B_{2m}$ are Bernoulli numbers; recall also that  
\begin{equation*}
\zeta( - n) =-\frac{B_{n+1}}{n+1}, 
\end{equation*}
 for odd $n = 1, 3, 5, \ldots$.
\begin{equation*}
\zeta( 1- 2m) = -\frac{B_{2m}}{2m}, \text{ if} \; m > 0.
\end{equation*}

\begin{example}
	(By Euler ),
\begin{equation*}
 \zeta(2) =  \frac{{\pi}^2}{6}. 
  \; \zeta(4) = \frac{{\pi}^4}{90}, \; \zeta(6) = \frac{{\pi}^6}{3^35\cdot7}, 
\end{equation*}
\begin{equation*}
\zeta(-1) = -\frac{B_{2}}{2} = -\frac{1}{12}, 
  \zeta(-3) = \frac{1}{120}. 
\end{equation*}
\end{example}

Define polylogarithm 
\begin{equation*}
L_m(z) = \sum_{n=1}^\infty z^n n^{-m}.
\end{equation*}
\begin{example}
\begin{equation*}
  \zeta(2) = L_2(1).
\end{equation*}
\end{example}

\section{ On class groups of rings with divisor theory on imaginary quadratic and cyclotomic fields}
The study of class groups  of  rings and corresponding schemes is an actual scientific problem 
(see~\cite{BorevichShafarevich1986,Hashimoto2017} and references therein).
For regular local rings, according to the Auslander-Buchsbaum theorem, the  (divisors) class group is trivial.
But in most interesting cases the group is nontrivial.
The Heegner approach, together with the results of Weber, Birch, Baker and Stark, makes it possible to calculate and even parametrize rings with a given (small) class number in some cases.
Let $R$ be a commutative ring with identity for which there exists the theory of divisors~ \cite{BorevichShafarevich1986}.
The order of the class group is calculated on the basis of the use of $ L $-functions.
We investigate one of the aspects of this problem, consisting in finding the moduli spaces of elliptic curves defined over the rings $R$ with the given class number.

{\bf Problem.}
To investigate the case of elliptic curves over  rings of integers of quadratic fields (rings of integers $\mathcal O$ of quadratic algebraic extensions $k$  of the field of rational numbers $\mathbb Q$) with a small  class number, see  
\cite{BorevichShafarevich1986}.

In some cases, for instance under computer algebra computations,  we have to enumerate investigated objects. Some simple parametric spaces and moduli spaces in the case of imaginary quadratic fields are presented below~\cite{Glazunov2018}.
We present an elementary introduction to this problem and give the moduli spaces as trivial bundles over affine part of the groups of rational points of some elliptic curves over the ring of integers $\mathbb Z$.
Below we  present parameter spaces and moduli for class number one and two.
Let 
\begin{center}
$E: y^2 = x^3 + ax + b, \; Disc(E) = 4a^3 + 27b^2, Disc(E) \ne 0, $                                \;                          (*) 
\end{center}
be an elliptic curve over the ring $\mathcal O$ .
Let $A_1$ be the affine part of the group of rational points over $\mathbb Z$ of the Heegner elliptic curve $y^2 = 2x(x^3 + 1)$.
With results by Heegner, Deuring, Birch, Baker, Stark, Kenku, Abrashkin, we deduce 
\begin{proposition}
Let $\mathcal O$ be the ring of integers of the imaginary quadratic field with class number one.
Then the parameter space of elliptic curves of the form (*) is the trivial bundle
\begin{center}
$({\mathcal O}\times{\mathcal O}/(Disc(E)=0))\times A_1$.
\end{center}
\end{proposition}
\begin{proposition}
Let $k$ be  the imaginary quadratic field with class number one.
Then the moduli space of elliptic curves of the form (*) is the trivial bundle
\begin{center}
$k\times A_1$.
\end{center}
\end{proposition}
Let $A_2$ be the affine part of the group of rational points over $\mathbb Z$ of the elliptic curve
$X^3 + 3Х = -Y^2$, let $A_3$ be the affine part of the group for the elliptic curve
$X^3 - 3Х = 2Y^2$, and $A_4$ respectively for
$9X^4 - 1 = 2Y^2$.
\begin{proposition}
Let $\mathcal O$ be the ring of integers of the imaginary quadratic field with class number two.
Then the parameter spaces of elliptic curves of the form (*), without  an exceptional case,  are trivial bundles
\begin{center}
$({\mathcal O}\times{\mathcal O}/(Disc(E)=0))\times A_2$,
$({\mathcal O}\times{\mathcal O}/(Disc(E)=0))\times A_3$,
$({\mathcal O}\times{\mathcal O}/(Disc(E)=0))\times A_4$.
\end{center}
\end{proposition}
\begin{proposition}
Let $k$ be  the imaginary quadratic field with class number two.
Then the moduli spaces of elliptic curves of the form (*), without  an exceptional case, are the trivial bundles
\begin{center}
$k\times A_2$, $k\times A_3$, $k\times A_4$.
\end{center}
\end{proposition}
\begin{theorem}
 (The Kronecker-Weber theorem)
Every finite abelian extension of $\mathbb Q$ is contained in a cyclotomic field.
\end{theorem}

With results by Heegner, Deuring, Birch, Baker, Stark, Shafarevich we have
\begin{proposition}
Imaginary quadratic fields with class number one and with descriminants $-D =  4, 8, 3, 7, 11, 19, 43, 67, 163$ are contained, respectively, in cyclotomic fields  
\begin{equation}
\begin{split}
\mathbb Q(\sqrt[ 4]{1}), \mathbb Q(\sqrt[ 8]{1}), \mathbb Q(\sqrt[ 3]{1}), 
\mathbb Q(\sqrt[ 7]{1}), \mathbb Q(\sqrt[ 11]{1}), \\
\mathbb Q(\sqrt[ 19]{1}), \mathbb Q(\sqrt[ 43]{1}), 
\mathbb Q(\sqrt[ 67]{1}),
\mathbb Q(\sqrt[ 163]{1}). 
\end{split}
\end{equation}
\end{proposition}

\section{Eisenstein Series}

Here we follow to~\cite{Weber1908,Hasse1931,Deuring,BorevichShafarevich1986}.

Let $\tau$ belong to the modular figure of the modular group $\Gamma =\Gamma(1)$.

\begin{definition}
In these notations with $k>1$ the Eisenstein series  is defined as
$$c_{k} = \sum_{m \ne 0, k > 1} \frac{1}{(n + m\tau)^{2k}}.$$
\end{definition}

\begin{proposition}
 Eisenstein series  have the representation

$$ c_{k} =  2  \zeta(2k) + \frac{2 (- 2\pi i)^{2k}}{(2k -1)!} \sum_{n> 0, m > 0} n^{2k-1} q ^{nm},$$

where $q = e^{2\pi i \tau} \ne 0$.

If we will use functions of the sums of divisors $\sigma_{2k-1}$ we obtain

$ c_{k} =  2  \zeta(2k) + \frac{2 (- 2\pi i )^{2k}}{(2k -1)!} \sum_{n=1}^{\infty}  \sigma_{2k-1}(n) q ^{n}$

 or shortly

$ c_{k} =  2  \zeta(2k) + \frac{2 (- 2\pi i)^{2k}}{(2k -1)!} S_{2k-1}$.
\end{proposition}
As $\zeta(2k) = (-1)^{k-1} \frac{(2\pi )^{2k}}{2(2k)!}B_{2k}$ we have
\begin{corollary}
$c_{k} = 2  \zeta(2k)(1 - \frac{4k}{B_{2k}}\sum_{n=1}^{\infty}  \sigma_{2k-1}(n) q ^{n}).$
\end{corollary}
Put $g_2 =  60 c_2$, $g_3 = 140 c_3$.

\begin{proposition}
 $\Delta = g_2^3  - 27g_3^2 \ne 0.$ 
\end{proposition}

As  $\Delta  \ne 0 $ it is possible to define $J = \frac{g_2^3}{\Delta}$.
\begin{definition}
Modular invariant of the elliptic curve $y^2 = 4x^3 - g_2 x - g_3$ is equal to $j = 2^63^3 J$.
\end{definition}
\begin{proposition}
$j = \frac{1}{q} + u_1 q + \cdots$
 where $u_i$ are integers, $u_0 = 0$.
\end{proposition}

Let us transform $ c_{k}$ in such a way that corresponding Fourier coefficients under 
$ q ^{n}, n \ge 1$ will rational numbers. Dividing  $ c_{k}$ on $2  \zeta(2k)$ and denoting the obtained result as $E_k$ we have by the Corollary 1 
\begin{proposition}
$E_k =   1 - \frac{4k}{B_{2k}}\sum_{n=1}^{\infty}  \sigma_{2k-1}(n) q ^{n}. $
\end{proposition}
\begin{example}
\begin{equation*}
E_2 =   1 + 240 \sum_{n=1}^{\infty}  \sigma_{3}(n) q ^{n} , 
\end{equation*}
\begin{equation*}
E_3 =   1 - 504\sum_{n=1}^{\infty}  \sigma_{5}(n) q ^{n}.
\end{equation*}
\end{example}

\section{Class group, class fields and zeta functions}
Here we follow to~\cite{Hasse1931,BorevichShafarevich1986}.

Let $K$ be an imaginary quadratic field and let $Cl_K$ be its class group.
\begin{definition}
Let $N(\mathfrak{a})$ be the norm of the ideal $\mathfrak{a}$.
 The Dedekind $\zeta$-function for $K$ is defined for all $s>1$ by the series \\
 $$\zeta_K(s) = \sum  \frac{1}{N(\mathfrak{a})^s },$$
where the sum is taken over all nonzero ideals $\mathfrak{a} \in \mathcal O_K$ .
\end{definition}

Let $R$ be a subring ($R \ne \mathbf Z)$  of the ring of integers $\mathcal O_K$ of the imaginary quadratic field $K$.

Let $M_1, \ldots M_h$ be pairwise nonequivalent modules of $K$ with the same ring of multipliers $R$. 

\begin{proposition}
 $j(M_1),\ldots , j(M_h)$  are integer algebraic numbers which are conjugate over $K$.
\end{proposition}

\begin{proposition}
The field $K(j(M_i))/K$ is the normal field. 
\end{proposition}
\begin{definition}
 The field $K(j(M_i))/K$ is called the ring class field.
\end{definition}

Follow to~\cite{Hasse1931} it is possible to define ray class field.
As in an  imaginary quadratic field there is no real infinite primes so modulus of the field is an ideal of the ring of integers of the field. 

Let $\mathfrak{m}$ be a modulus of the an  imaginary quadratic field $K$, let $Cl_K^{\mathfrak{m}}$ be the ray class group, let  
$\tau_W$ be the Weber function .

Let $\mathfrak{R} \in Cl_K$ and let $\mathfrak{R^*} \in Cl_K^{\mathfrak{m}}$ be the  ideal class whose image in $Cl_K$ is equal to
$({\mathfrak{m}})\mathfrak{R}^{-1}$.
\begin{proposition}
The field $K(j(\mathfrak{R}), \tau_W(\mathfrak{R^*}))/K$ is the ray class field. 
\end{proposition}

Let $C$ be an ideal class.

\begin{definition}
The	ideal class zeta function is  the expression of the form
$$\zeta_C(s) = \sum_{\substack{\mathfrak{a} \in C\\
\mathfrak{a} \; \text{integral}}} \frac{1}{N(\mathfrak{a})^s } $$
$$ s = \sigma + it, \; \sigma > 1.  $$  
\end{definition}

   Below we present values of zeta and $L$-functions connecting with  imaginary quadratic fields.
Let $d$ be a squarefree integer number, $ K = \mathbb Q(\sqrt{d}) $  a quadratic field, $\chi$ be the character of the quadratic field $K$.
Let $L(s,\chi)$ be the $L-$series with a nonunit character $\chi$ modulo $|D|$. Here $D$ is the discriminant of the field $K$.
\begin{proposition}
\begin{equation*}
 \zeta_K(s) = \zeta(s) L(s,\chi) =  \zeta(s) \prod_p (1 -  \frac{\chi(p)}{p^s }).
\end{equation*}
\end{proposition}

Let $m$ be the number of roots of unity of the imaginary quadratic field $[K : \mathbb Q]$. 
\begin{remark}
 $m = 4 $ for $ K = \mathbb Q(\sqrt{-1}),  \; m = 6 $ for $ K = \mathbb Q(\sqrt{-3}), \; m = 2$ for all other imaginary quadratic fields.
\end{remark}
Let $h$ be the class number of the field $K.$
\begin{proposition}
\begin{equation*}
  L(1,\chi) = \frac{2\pi h}{m \sqrt{|D|} }.
\end{equation*}
\end{proposition}
\begin{corollary}
For  imaginary quadratic fields with class number one ($h = 1$) we have
$$
  L(1,\chi) = \left\{
   \begin{array}{lll}
      \frac{\pi}{4}, K = \mathbb Q(\sqrt{-1})\\
       \frac{\pi}{3 \sqrt{3}},  K = \mathbb Q(\sqrt{-3})\\
        \frac{\pi}{\sqrt{|D|}}, -D = 8,7,11,19, 43, 67, 163. \\
    \end{array}
     \right.
$$
\end{corollary}

\section{Multiple zeta  values} 

\begin{definition}
  Let $x_1, \ldots x_p$ be natural numbers with $x_p \ge 2$.
The multiple zeta value of the weight $w$ and the depth $p$ is called the expression of the form
\begin{equation*}
\zeta(x_1, \ldots x_p) = \sum_{0 < n_1 < \cdots < n_r}  \frac{1}{n_1^{x_1}\cdots n_p^{x_p}}, w= \sum x_i .
\end{equation*}
\end{definition}

\begin{example}
\begin{equation*}
	\zeta(2, 2) = \sum_{0 < n_1 <  n_2}  \frac{1}{n_1^{2} n_2^{2}} , w= \sum x_i = 4.
\end{equation*}
\end{example}

\begin{example}
\begin{equation*}
	\zeta(2, 2) =  \frac{1}{2} (\zeta(2)\zeta(2) - \zeta(4)). 
\end{equation*}
\end{example}

Let $\mu_N$ be the group of roots of unity.
\begin{definition}
 Let $x_1, \ldots x_p$ be natural numbers with $x_p \ge 2$.
The multiple zeta value relative to $\mu_N$  of the weight $w$ and the depth $p$ is called the expression of the form \\
\begin{equation*}
\zeta(\substack{x_1, \ldots x_p \\
\epsilon_1, \ldots, \epsilon_p }  ) = \sum_{0 < n_1 < \cdots < n_p}  \frac{\epsilon_1^{n_1} \ldots, \epsilon_p^{n_p}}{n_1^{x_1}\cdots n_p^{x_p}} ,  \epsilon_i \in \mu_N,
\end{equation*}
\begin{equation*}
  w= \sum x_i ,  
(x_p, \epsilon_p) \ne (1,1).
\end{equation*}
\end{definition}


\section{Elements of   non-Archimedean local fields and $ p-$adic  analysis}

Here we  present elements of $p-$adic local fields, their algebraic extensions and $p-$adic interval analysis. We follow to \cite{BorevichShafarevich1986,Koblitz}.

\subsection{Elements of   non-Archimedean local fields}
A  non-Archimedean local field is a complete discrete valuation field with finite residue field. Further, for brevity, we call these fields local.
In other words, a field $K$ is called local if it is complete in a topology determined by the valuation of the field and if its residue field $k$ is finite. We assume further that the valuation ${\nu}$ is normalized, i.e. the homomorphism of the multiplicative group of the field to the additive group of rational integers ${\nu}: K^*  \to {\mathbb Z}$ is surjective.

The structure of such fields is known: if the field $K$ has the characteristic zero, then it is a finite extension of the $p-$adic field 
${\mathbb Q}_p$, which is the completion of the field of rational numbers with respect to the $p-$adic valuation.

If $[K:{\mathbb Q}_p] = n$, then $n = ef$, where 
$f$ is the degree of classes of residues, (i.e. $f = [k:{\mathbb F}_p]$) and 
$e = \nu_K(p)$ is the ramification index of  $K$..

If the field $K$ has the characteristic $p$, then it is isomorphic to the field $k((T))$ of formal power series, where $T$ is  a uniformizing parameter.

Let $L$ be a finite extension of a local field $K$ with their residue fields $l$ and $k$, $p = char \; k$  and $e_{L/K}$ be the ramification index of $L$ over $K$.

An extension $L/K$ is called unramified if a) $e_{L/K}= 1$; b) the extension $ l/k $ is separable. An extension $L/K$ is called tamely ramified if a) $p$  does not divide $e_{L/K}$; b) the extension $ l/k $ is separable.

An extension $L/K$ is called  wildly ramified if  $e_{L/K} = p^s, \; s \ge 1$;

Denote  by $Tr_{L/K}$ and by $Norm_{L/K}$ respectively the trace and the norm of the extension ${L/K}$. We drop indices, when it is clear what kind of extension we are talking about.

Denote by $K_{nr}$ the maximal unramified extension of the field $K$ (in a fixed algebraic closure of the field $K$) with a residue field $k_s$, which is the algebraic closure of a field $k$.

In a non-Archimedean local field $K$ each of its elements $\alpha$ has a representation 
$\alpha = \epsilon \pi^m$, where $\epsilon$ is a unit of the ring of integers of the field $K$ and $\pi$ its uniformizing element, that is ${\nu}(\pi) = 1$  , $m$ is an integer rational number.
A unit is called principal if $\epsilon \equiv 1\pmod \pi $.

\begin{lemma}
If the local field contains a primitive $p-$th root $\xi_p$ of unity, then $\nu(\xi_p - 1) = \frac{e}{p-1}$ is an integer number.
\end{lemma}
{\bf Proof}. $\xi_p - 1$  is the root of the equation 
$(x+1)^{p-1} + (x + 1)^{p-2} + \cdots + (x+1) + 1 = x^{p-1} + p(\cdots) + p.$ The value of the $p-$adic valuation at the root of this equation is $\frac{e}{p-1}$  which proves the required. $\Box$

A complete discrete valuation field with an algebraically closed residue field is called a quas-ilocal field.

\subsection{$p-$adic intervals and $p-$adic distributions}

Let $X$ be a topological space. A distribution on $X$ with values in an abelian group $A$ is a finitely additive function from the compact-open subsets of $X$ to $A$. Let $|\;|_p$ be the $p-$adic norm.

Define $[\alpha, N]_p = \{x \in {\mathbb Q}_p | |x - \alpha|_p \le \frac{1}{p^N}\}$,
$ \alpha \in  {\mathbb Q}_p, N \in {\mathbb N}$.

\begin{definition}
We call sets $[\alpha, N]_p$ the $p-$adic intervals (disks) and define by these $p-$adic intervals the basis of open sets 
on ${\mathbb Q}_p$.
\end{definition}
It is easy to test that axioms of open sets are satisfied.
\begin{remark}
$p-$adic intervals $[\alpha, N]_p$ open and closed simultaneously.
\end{remark}
{\bf Proof}. Any union of open $p-$adic intervals is open. Intervals $[\alpha, N]_p$ are closed, because $[\alpha, N]_p$ is an addition to the union of open intervals $[\alpha^`, N]_p$ for all $ \alpha^` \in  {\mathbb Q}_p$ for which $\alpha^` \neg \in [\alpha, N]_p$. $\Box$

Further we will call $[\alpha, N]_p$ as intervals.  More generally we will consider compact-open sets. Let $X$ be a compact-open set. Recall that a function $f: X \to {\mathbb Q}_p$ is is locally constant if and only if $f$ has a representation as a finite linear combination of characteristic functions of compact-open subsets.

Let $U = U_1 \cup U_2 \cup \cdots \cup U_n $ be a partition of $U \subset X$.
Recall that the additive mapping $\mu$ of a set of  compact-open subsets of $X$ with value in ${\mathbb Q}_p$  is called the 
$p-$adic distribution on $X$: 
$$\mu(U) = \mu(U_1) + \mu(U_2) + \cdots + \mu(U_n) .$$

\subsubsection{Bernoulli  distributions.} Let $B_m(x)$ be the $m-$Bernoulli polynomial. These polynomials are defined by the decomposition
$$\frac{te^{xt}}{e^t - 1} = \sum_{m=0}^{\infty} B_m(x) \frac{t^m}{m!}.$$ 
We have: $B_0(x) = 1, \; B_1(x) = x - \frac{1}{2}, \; B_2(x) = x^2 -x +\frac{1}{6}, \; 
B_3(x) = x^3 - \frac{3}{2}x^2 + \frac{1}{2}x, \; B_4(x) = x^4 - 2x^3 +x^2 - \frac{1}{30}, \ldots $

\begin{remark}
If we substitute $x = 0$ in the $m-$Bernoulli polynomial  we obtain $m-$Bernoulli number:
\begin{equation*}
 B_0(0) = 1, \; B_1(0) =  - \frac{1}{2}, \; B_2(0) = \frac{1}{6}, \; 
 B_3(0) = 0, \; B_4(0) =  - \frac{1}{30}, \ldots 
\end{equation*}
\end{remark}

Let now for $\alpha$ the inequality $0 \le \alpha \le p^N - 1$ is satisfied.
Define the function $\mu_{B,m}$ by the formula
$$\mu_{B,m}([\alpha, N]_p) = p^{N(m - 1)}B_m(\alpha / p^N)  .$$ 

\begin{proposition}
The function $\mu_{B,m}$  is expanded to the distribution on ${\mathbb Z}_p$. This distribution for the given $m$ is called the 
$m-$th Bernoulli  distribution.
\end{proposition}

\newpage

\section{$ p-$adic  interpolation of zeta and $L$-functions}
Here we follow to \cite{Leopoldt1959,Leopoldt1960,Leopoldt1962,KubotaLeopoldt1964,Leopoldt1975,Iwasawa2001,BorevichShafarevich1986}
but consider  Leopoldt's zeta functions and Leopoldt's $L$-functions only. 
\subsection{Leopoldt's zeta functions}
Recall that  for natural $n > 0 $ and for Bernoulli numbers in pair notation
\begin{equation*}
\zeta( 1- n) = -\frac{B_{n}}{n}.
\end{equation*}
Recall the generalization of the Kummer congruences.
Let $ p$ be a prime number.
\begin{lemma}
Let $m, n$ be pair natural numbers such that $p$ not devides $m, $ $p$ not devides $n, $ and 
$m \equiv n  \bmod (p - 1) p^{N }.$ Then
\begin{equation*}
(1 - p^{m - 1}) \frac{B_{m}}{m} \equiv (1 - p^{n - 1})\frac{B_{n}}{n} \bmod p^{N + 1}.
 \end{equation*}
\end{lemma}
Let $p \ge 3$ be a prime number and $P_p$ the set of pair natural numbers $0, 2, \ldots, p - 3.$
For a given $P_p = \{ 0, 2, \ldots, p - 3\}$ and $a \in P_p$ denote by 
$C_a^p$ 
the set of natural numbers $m$ such that 
 $ m  \equiv a \bmod (p-1)$ .
\begin{lemma}
Each $C^p_a$ is dense in the set of $p-$adic integer numbers ${\mathbb Z}_p$.
\end{lemma}

By Leopoldt and others this gives possibility to expand  the function
\begin{equation*}
 -(1 - p^{m - 1})\zeta(1 - m)
 \end{equation*}
on the ring of $p-$adic integer numbers.

Indeed it is possible to expand the function $ -(1 - p^{-s})\zeta(s)$ which is defined at points $1 - m, m \in C^p_a$ on all 
 ring ${\mathbb Z}_p$.
Denote these  functions by $\zeta_{p,a}(s), \; s \in {\mathbb Z}_p.$ 
We will call this expansion the $p-$adic continuation.
\begin{proposition}
 Let $ -(1 - p^{-s})\zeta(s)$ be zeta with Euler multiplier
 
$(1 - p^{-s}).$ There are $\frac{p - 1}{2}$ $p-$adic continuations $\zeta_{p,a}(s)$ of the function

 $ -(1 - p^{-s})\zeta(s). $ For $a = 2, \ldots, p - 3 $ these continuations are $p-$adic analytic functions and for $a = 0$
 this continuation $\zeta_{p,0}(s)$ is $p-$adic meromorphic function with a single pole of the 1st order at the point $s = 1.$
\end{proposition}


\subsection{Leopoldt's $L$-functions}
Let $\chi$ be the numeric character modulo $m > 1$ and $n \ge 1$ natural number. Then
\begin{equation*}
L( 1- n,\chi) = -\frac{B_{n,\chi}}{n}.
\end{equation*}
where $L$ is the Dirichlet $L-$function and  ${B_{n,\chi}}$ are generalized Bernoulli numbers (1).
If we will use Bernoulli polynomial $B_n(x),$ then
\begin{equation*}
B_{n,\chi} = m^{n - 1} \sum_{r=1}^m \chi(r) B_n(\frac{r}{m}).
\end{equation*}
\begin{remark}
For any integer rational $c, \; (c,p) = 1$ the sequence $\{c^{n^p}\}$ converges in the field ${\mathbb Q}_p$.
Let $\gamma = \lim_{n \to\infty} c^{n^p}$. Then $\gamma \equiv c \bmod (p)$ , $\gamma^{p - 1} = 1.$
\end{remark}
\begin{proposition}
For any primitive numeric character $\chi$ and any prime $p$ there exist $p-$adic function $ L_p(1 - n, \chi)$ defined at 
integer $p-$adic numbers $s$ (with the exaption $s = 1$ in the case of unit character $\chi = 1$ ) with the property
\begin{equation*}
 L_p(1 - n, \chi) = -(1 - \chi(p) \gamma^{ - n}(p)\cdot p^{n - 1})\frac {B_{n,\chi \gamma^{ - n}}}{n}, \; n \ge 1.
\end{equation*}
For odd character $\chi$  the function $ L_p(1 - n, \chi)$ equals to identically zero.

\end{proposition}
\subsection{Kubota-Leopoldt L-functions}

\subsection{Tate module of a number field}

\subsection{On Iwasawa conjecture}

\subsection{Constructions of p-adic $L-$functions}

\newpage
\section{ Iterated integrals and (multiple) zeta values}
Here we follow to\cite{Parshin,Chen}.

Let ${\mathbb C}$ be the complex plane and $f_i(z) $ be the holomorphic function on ${\mathbb C}$ . Let $f_i(z)dz $  be the differential of the first kind on ${\mathbb C}$. Let $S$ be a Riemann surfaces and $w$ be the differential of the
 first kind on $S$.
Parshin  has considered iterated integrals of this type on Riemann surfaces~\cite{Parshin}.
Chen~\cite{Chen} for  smooth paths on a manifold $M$ and respective path spaces have investegated  iterated (path) integrals. For differential forms $w_1, \ldots , w_r$ on $M$ he has constructed the iterated integrals by repeating $r$ times the 
integration of the path space differential forms  (and their linear combinations). Chen~\cite{Chen} has denoted the iterated integrals as
$\int w_1 w_2 \cdots w_r$ and set $\int w_1 w_2 \cdots w_r = 1$ when $r = 0$ and  $\int w_1 w_2 \cdots w_r = 0$ when $r < 0$.

\begin{example}
\begin{equation*}
\zeta(2) = \int_0^1 \frac{dt_1}{t_1} \int_0^{t1} \frac{dt_2}{1 - t_2} = \frac{{\pi}^2}{6}.
\end{equation*}
\end{example}
More generally iterated integrals are path space differential forms which permit further integration.

\section{Formal groups and $p$-divisible groups}
Recall some definitions.
Let $K$ be a complete discrete variation field with the ring of integers $O_K$ and the maximal ideal $ M_K$. 
A complete discrete variation field with finite residue field is called a {\it local} field~\cite{CF:ANT}. 
A complete discrete variation field $K$ with algebraically closed residue field $k$ is called a {\it quasi-local} field~\cite{Vv:QEC}. Below we will suppose that in the case the characteristic of $k$ satisfies $p > 0$.
Let $K$ be a local  or quasi-local field. If $K$ is a local field~\cite{CF:ANT} and  has the characteristic $0$ then 
it is a finite extension of the field of $p$-adic numbers ${\mathbb Q}_p$. Let ${\nu}_K$ be the normalized exponential valuation of $K$. If $[K : {\mathbb Q}_p ] = n$ then $n = e \cdot f$, where $e = {\nu}_{K}(p)$ and $f = [k : {\mathbb F}_p ]$, where $k$ is the residue field of $K$ (always assumed perfect ). If $K$ has the characteristic $p > 0$ then it  isomorphic to the field $k((T))$ of formal power series, where $T$ is uniformizing parameter. 
Let $L$ be a finite extension of a local field $K$, $k, l$ their residue fields, $p = char\; k$ and $e_{L/K}$ ramification index of $L$ over $K$. 
An extension $L/K$ is said to be ${\it  unramified}$ if $e_{L/K} = 1$ and extension $l/k$ is separable.
An extension $L/K$ is said to be ${\it tamely \;  ramified}$ if $p $ not devides $e_{L/K}$  and the residue extension $l/k$ is separable.
An extension $L/K$ is said to be ${\it totally \;  ramified}$ if $e_{L/K} = [L:K] = (char\; k)^s, s \ge 1$.

Let $L/K$ be the finite Galois  extension of quasi-local field $K$ with Galois group $ G$, $F(x,y)$ one dimensional formal group low over the ring of integers $O_K$ of the field $K$, $F( M_K)$ be the $G$ - module, that is defined by the group low $F(x,y)$ on the maxilal ideal $ M_K$ of the ring $O_K$, $\ M_{K}^t (t \in {\mathbb Z}, t \ge 1)$ be the subgroup of $t$-th degrees of elements from $ M_K, F^t_K := F({ M}_K^t)$.

\begin{definition}
\label{Mu}
For $n \in {\mathbb Z}$ the function $\mu (n)$,  $N_{L/K} (F^{n}_L) \subset F^{\mu (n)}_K$ 
is defined by the condition: $F^{\mu (n)}_K$ is the least of subgroups $F^t_K $ ($t = 1, 2, \ldots$) containes $N_{L/K} (F^n_L)$.
\end{definition}
\begin{remark}
Please do not confuse with the measure $\mu$.
\end{remark}
Below we will suppose that $char\; k > 3$.

\subsection{Norm Maps}

Here we use results on formal groups from    \cite{Vv:UNFG, Ha:FG}.
Let $F_L =  F( M_L)$ be the $G$ - module that is defined by the $n$-dimentional group low $F(x,y)$  on the product 
$ (M_L)^n := M_L\times \cdots \times M_L$,  ($n$ times) of maximal ideals of the ring $O_L$ of any finite Galois extension $L$ of the field $K$.

\begin{definition}
The norm map
$N: F_L \to F_K$
of the module $F_L$ to $F_K$ is defined by the formula
$N(a) = (((a +_{F} \sigma a) +_{F} \cdots) +_{F} \sigma_s a)$,
where $a +_{F} b$ denotes the addition of points in the sense of group structure of the module $F_L$,
$a, b \in M_L$, $G = Gal(L/K)$, $\sigma_s \in G$, $[G : 1] = s$.
\end{definition}

Let $p := char \; k, \; e := {\nu}_K (p)$, ($e = +\infty$, if characteristic  of the field $K$ is equal $p$ and $e$ is positive integer in the opposide case), $L/K$ be the Galois extension of the prime degree $q$, $F(x,y)$ be the one dimensional group low over $O_K$. Let $p := char \; k > 0$.

\begin{lemma}
  If $\Pi_{s} \in \pi_L^s \cdot O_L $, $s \ge 1$ then

$N(\Pi_{s}) \equiv Tr(\Pi_{s}) + \sum_{n = 1}^{\infty} c_{n} [Norm \; \Pi_{s}]^n (mod \; Tr (\pi_{L}^{2s} \cdot O_L))  $

where $c_n \in O_K$ are coefficients of the $p$ - iteration of the group low. 
\end{lemma}

  Let $R$ be a commutative ring. 
Let $A, \; B, \; C$ be finite group schemes over $R.$
The sequence
$$
\begin{array}{ccccccc}
0&\longrightarrow & A &
\stackrel{f}{\longrightarrow}& B &
\stackrel{g}{\longrightarrow} & C
\end{array}
$$
is called exact if $Im \; f = Ker \; g.$

Let $p$ be a prime number and $h$ be an integer, $h \ge 0.$

Recall the definition of the $p$-divisible group by J. Tate.
\begin{definition}
A $p$-divisible group over $R$ of height $h$  is an inductive system
\begin{equation*}
  G = (G_{\nu}, i_{\nu}), {\nu} \ge 0,
\end{equation*}
where \\
(i)  $G_{\nu}$ is a finite group scheme over $R$ of order $ p^{\nu h}$,\\
(ii) for each ${\nu} \ge 0,$
$$
\begin{array}{ccccccc}
0&\longrightarrow & G_{\nu} &
\stackrel{i_{\nu}}{\longrightarrow}& G_{\nu + 1} &
\stackrel{p^{\nu}}{\longrightarrow} & G_{\nu + 1}
\end{array}
$$
is exact.
\end{definition}
  
\begin{example}
  $(G_{\nu} = ({\mathbb Z}/p^{\nu}{\mathbb Z})  )^h, \; G = ({\mathbb Q}_p /{\mathbb Z}_p)^h.$   
\end{example}

\section{Motives and  ($p$-adic) (multiple) zeta values}

Glanois in paper~\cite{Glanois2016} presents the revised and expanded version of his Doctoral thesis  
 [Periods of the motivic fundamental groupoid of $ \mathbb{P}^1 
\{0,\mu_N, \infty \}$, Pierre and Marie Curie University, 2016;], written under F. Brown.

Let $k_N = \mathbb{Q}(\xi_N)$ be the cyclotomic field, $\xi_N \in \mu_N$ be a primitive $N$th root of unity and ${\mathcal O}_N$ be the ring of integers of $k_N$.
The corresponding multiple zeta values at arguments $x_i \in \mathbb{N}, \epsilon_i \in \mu_N$ can be expressed in terms of the coefficients of a version of Drinfeld`s assosiators by Drinfeld~\cite{Drinfeld:DA}, which in turn, can be expressed in terms of periods of the corresponding motivic multiple zeta values (MMZV).

These MMZV $\zeta^{\mathfrak m}(\substack{x_1, \ldots x_p \\
\epsilon_1, \ldots, \epsilon_p }  ) , \;
 \epsilon_i  \in \mu_N, \; (x_p, \epsilon_p) \ne (1,1)$ 
relative to $ \mu_N$   (of the weight $w = \sum x_i $ and the depth $p$), are elements of an algebra $\mathcal H^N$ over $\mathbb{Q}$ and span the algebra.

The algebra $\mathcal H^N$ carries an action of the motivic Galois group of the category of mixed Tate motives over 
${\mathcal O}_N[1/N]$.
The author studies the Galois action on the motivic unipotent fundumental groupoid of $ \mathbb{P}^1 \backslash \{0,\mu_N, \infty \}$ (or of  $\mathbb{G}_m \backslash \mu_N$) for next values of $N: N \in \{ 2^a3^b, a + 2b \le 3\} = \{1, 2, 3, 4, `6`, 8 \}  $.

His results include: bases of multiple zeta values via multiple zeta values at roots of unity $\mu_N$ for the above $N$;  more generally, constructing of families of motivic iterated integrals with prescribed properties;
the new proof, via the coproduct by Goncharov~\cite{Goncharov2005} and its extension by Brown~\cite{Brown2012}, of the results by Deligne~\cite{Deligne2010} that the Tannakian category of mixed Tate motives over ${\mathcal O}_N[1/N]$ `for $N = \{2, 3, 4, 8\}$ is spanned by the motivic fundumental groupoid of $\mathbb{P}^1 \backslash \{0, \mu_N,  \infty \}$ with an explicit basis`.  

In article \cite{Unver2016} Unver continues his investigation of $p$-adic multiple zeta values\cite{Unver2004},
 presenting a computation of values of the p-adic multiple polylogarithms at roots of unity.
  The main result of the paper \cite{Unver2016} (Theorem 6.4.3 with Propositions 6.4.1 and 6.3.1) is to give explicit expression for the cyclotomic $p$-adic multi-zeta values $\zeta_{p}(s_1, s_2; i_1, i_2)$ of depth two.
The result is far too technical to state here.

The proof of the theorem is rather technical; it is based on rigid analytic function arguments and a long distance analysis of group-like elements of related algebras. 

For number fields the category of realizations has defined and investigated by Deligne~\cite{Deligne1989}.
Results  on multiple zeta values, Galois groups and geometry of modular varieties has presented by 
Goncharov~\cite{Goncharov2000}.
Interesting unipotent motivic fundamental group is defined and investigated by Deligne and
 Goncharov~\cite{DeligneGoncharov2005}.
Tannakian interpretation of $p$-adic multiple zeta values is given by Furusho~\cite{Furusho2007}.

Results obtained in the paper \cite{Unver2016} may be applied to the problems of the $p$-adic theory of higher cyclotomy.

\section{On the Eisenstein series associated with Shimura varieties}
  Interesting classes of Shimura varieties form varieties which have an interpretation as moduli spaces of abelian varieties.
Moduli spaces of corresponding $p-$divisible groups over perfect fields of characteristic $p$ are used for investigation of the local structure of such Shimura varieties.

\subsection{On some Shimura varieties and their local structure}

Let $\mathbb{C}^n$ be $n-$dimensional complex vector space, $\{e_1, \ldots, e_{2n} \}$ be $2n-$linear independent vectors and
\begin{equation*} 
\Lambda = \{e_1 \mathbb{Z}, \ldots, e_{2n} \mathbb{Z} \}
\end{equation*}
 be a lattice.
Then 
\begin{equation*} 
\mathbb{C}^n / \Lambda
\end{equation*}
 is a compact commutative topological group.
If
 $\alpha \in GL_n(\mathbb{C})$
 and 
$\Lambda_1 = \alpha \Lambda$ 
then
\begin{equation*} 
 \mathbb{C}^n / \Lambda = \mathbb{C}^n / \Lambda_1.
\end{equation*}
If $n > 1$ then not for every lattice $\Lambda$ there exists an abelian variety.

\begin{proposition}
Let $ \mathbf{x}, \mathbf{y} \in \mathbb{C}^n$ and let $F(\mathbf{x}, \mathbf{y})$ be $\mathbb{R}-$bilinear form such that
\begin{equation*} 
(i)   F(\mathbf{x}, \mathbf{y}) = - F(\mathbf{y} , \mathbf{x}),
\end{equation*}
(ii) $ F(\mathbf{x}, i \mathbf{x})$ is the Hermitian positive defined form \\ 
and for \\
 (iii) $\mathbf{x}, \mathbf{y} \in \Lambda$ 
 it takes integer values:  $F(\Lambda,\Lambda) \in \mathbb{Z}$. \\
Then for this lattice $ \Lambda$ there exists the abelian variety.
\end{proposition}

\begin{definition}
The pair $(\Lambda, F)$ is  called the polarized abelian variety.
\end{definition}

Let $M = M(F)$ be the matrix of the form $F.$
\begin{definition}
The abelian variety $A = (\Lambda, F)$ is called principally polarized if the bilinear form
$F$ is unimodular or, equivalently, $det(M) = det(M(F)) = 1.$
\end{definition}

Denote by $\Pi$ the period matrix of the abelian variety $A.$ This is $n \times 2n$ complex matrix $\Pi = (M_1,M_2)$ with nondegenerate $n \times n$ matrices $M_1$ and $M_2.$

\begin{definition}
The period matrix $\Pi$ is called normalized if it has the form $(E_n,Z)$ where $E_n$ is the unit $n \times n$ matrix and
$Z \in \mathbb{H}_n,$  where.
\begin{equation*} 
\mathbb{H}_n = \{  n \times n \; \text{matrices} \; Z| Z^T = Z \; \text{and} \; Im Z > 0\},
\end{equation*}
is the Siegel upper half-plane. Here $Z^T$  is the  matrix transposed to $Z.$
\end{definition}

\begin{remark}
It is clear that the  Siegel matrix $Z \in \mathbb{H}_n$ defines the normalized period matrix $\Pi.$
\end{remark}
  Let 
$$
  J =   \left(
   \begin{array}{cc}
     0 & E_n\\
    - E_n & 0
           \end{array}
     \right).
$$
\begin{definition}
 Siegel  modular group $\Gamma_n = Sp_n(\mathbb{Z})$ is the 
set of matrices $$
  M = \left(
   \begin{array}{cc}
     A & B\\
     C & D
           \end{array}
     \right).
$$
such that 
\begin{equation*} 
  M^T J M = J.
\end{equation*}
\end{definition}

\begin{definition}
 Siegel  modular group $\Gamma_n$ acts on the the Siegel upper half-plane $\mathbb{H}_n $ by the formula
\begin{equation*} 
Z \mapsto (AZ + B)(CZ + D)^{-1}.
\end{equation*}
\end{definition}

\begin{proposition}
  In the framework of Definitions 13, 14 two Siegel matrices define isomorphic principally polarized abelian varieties if and only if one of them can be obtained from the other by the transformation of the Definition 15.
\end{proposition}

Sometimes we will use for Siegel matrices $Z$ an equivalent notations: $Z = X + iY, \; X, \; Y$ are real matrices, 
$X^T = X, \; Y^T = Y, \; Y > 0$; $Y$ is the matrix of the positive definite quadratic form.

\begin{definition}
Let $f$ be an analytic function on the Siegel upper half-plane that satisfies the equality
\begin{equation*} 
 f((AZ + B)/(CZ + D)) = (CZ + D)^h f(Z),
\end{equation*}
and is bounded on the domains of the form
\begin{equation*} 
\left\{ 
Z = X + iY, \; Z \in \mathbb{H}_n, Y \ge c E_n, \;  c > 0.
\right\}
\end{equation*}
Then the function $f$ is called Siegel modular form of the genus $n$ and the weight $k.$
\end{definition}

\begin{remark}
 As the matrix
$$
  J =   \left(
   \begin{array}{cc}
      E_n & B\\
    0 & E_n
           \end{array}
     \right).
$$
belong to $\Gamma_n$ (its determinant is equal 1) and
\begin{equation*} 
f(Z + B) = f(Z),
\end{equation*}
so $f(Z)$ has a representation by the Fourier series.
\end{remark}

\subsection{ On improper intersections of Kudla-Rapoport divisors and Eisenstein series by Sankaran}

Let $k$ be an imaginary quadratic field, $o_k$ its ring of integers and $o_{k,p}$ be the ring of integers of the completion $k_p$ of $k$ at $ p$. Sankaran~\cite{Sankaran2017} proves that the arithmetic degrees of Kudla-Rapoport cycles on an integral model of a Shimura variety attached to a unitary group of signature $(1,1)$ are Fourier coefficients of the central derivative of an Eisenstein series of genus 2. The main results of the paper are the following Theorem 4.13 on the value of the Eisenstein series and the Corollary 4.15 on the relation between  the arithmetic degree of special cycle and the Eisenstein series. These results confirm conjectures by Kudla~\cite{Kudla} and by Kudla, Rapoport~\cite{KudlaRapoport}
 on relations between intersection numbers of special cycles and the Fourier coefficients of automorphic forms in the degenerate setting and for dimension 2. As have pointed out by Kudla~\cite{Kudla2004}
and others `these relations may be viewed as an arithmetic version of the classical Siegel-Weil formula, which identifies the Fourier coefficients of values of Siegel-Eisenstein series with representation numbers of quadratic forms`. In the paper by Kudla, Rapoport 
~\cite{KudlaRapoport2010}  
`the Shimura variety is replaced by a formal moduli space of $p$-divisible groups, the special arithmetic cycles are replaced by formal subvarieties, and the special values of the derivative of the Eisenstein series are replaced by the derivatives of representation densities of hermitian forms.`   Sankaran defines the local Kudla-Rapoport cycles $Z(b)$ and gives some applications of results obtained in his earlier paper~\cite{ Sankaran2013} where he proved the Theorem 3.14 on cycles $Z(b)$. He allows `the polarizations to be non-principal in a controlled way`. An unpolarized case of $p$-divisible groups with the given $p^m$-kernel type and with applications to their Newton polygons has considered in the paper by Harashita~\cite{Harashita}.  
  Sankaran`s paper~\cite{Sankaran2017} consists of four sections. The first section presents the purpose of the paper and short description of ideas and results of next sections. Second section concerns with local Kudla-Rapoport cycles on the Drinfeld upper half-plane. The main result of this section is the Theorem 2.14 on values of local intersection numbers of these cycles. The third section is devoted to the prove of the closed-form formula for representation densities $\alpha(S,T)$. Author specializes the explicit formula on Hermitian representation densities by Hironaka~\cite{Hironaka}
to the case at hand: $F(S,T,X) \in Q[X], T \in Herm_2(o_{k,p}), ord_p det(T)$ is even, $ S = diag(p,1), \alpha(S_r,T) = F(S,T,(-p)^{-r})$. In the last section global aspects are discussed and main result is presented. Let $M^d_{(1,1)}$ denote the Deligne-Mumford (DM) stack over $o_k$ of almost-principally polarized abelian surfaces and ${\mathcal E}$ the DM stack over $o_k$ of principally polarized elliptic curves with multiplication by $o_k$. In conditions of the subsection 4.1 of the paper~\cite{Sankaran2017} author sets ${\mathcal M } = {\mathcal E} \times_{Spec o_k} M^d_{(1,1)}$ and define for $T \in Herm_2(o_{k})$ cycles ${\mathcal Z}(T)$.  Then in subsection 4.3 the author prove Theorem 4.13 and Corollary 4.15.  


\section*{Conclusions}
Classical and novel results on ($p$-adic) $L$-functions, ($p$-adic) (multiple) zeta values and Eisenstein series are presented.



\newpage
\end{document}